\documentclass[a4paper,fleqn]{cas-sc}

\usepackage[numbers]{natbib}

\def\tsc#1{\csdef{#1}{\textsc{\lowercase{#1}}\xspace}}
\tsc{WGM}
\tsc{QE}
\tsc{EP}
\tsc{PMS}
\tsc{BEC}
\tsc{DE}

\begin{document}
\let\WriteBookmarks\relax
\def\floatpagepagefraction{1}
\def\textpagefraction{.001}
\shorttitle{A criterion for grid deformation}
\shortauthors{Y. Zhen et~al.}

\title [mode = title]{An efficient and rigorous criterion for the admissible time step of deforming grids in generalized Lagrangian and ALE methods}                      



\author[1]{Yicun Zhen}
\cormark[1]
\ead{zhenyicun@proton.me}

\credit{Conceptualization, mathematical derivation, writing}

\affiliation[1]{organization={College of Oceanography, Hohai University},
                city={Nanjing},
                postcode={210024}, 
                state={Jiangsu},
                country={China}}

\author[1]{Shouxian Zhu}
\credit{Conceptualization, discussion, writing}

\author[2]{Shipeng Mao}

\credit{Discussion, suggestion of Prokhorova's theory, writing}

\affiliation[2]{organization={SKLMS, Academy of Mathematics and Systems Science, Chinese Academy of Science},
                city={Beijing},
                postcode={100190},
                country={China}}

\cortext[cor1]{Corresponding author}

\begin{abstract}
It is mathematically proved that, for a grid under a boundary-preserving 
 piecewise-smooth deformation, the positivity of the Jacobian determinant in each closed element is also a sufficient condition for the deformed grid to remain valid. Most of the proof is devoted to showing that the element-wise positivity of the Jacobian indeed induces a local isomorphism in the grid-deformation setting. This result provides users of ALE methods and of the generalized Lagrangian method with a clean answer: checking the Jacobian alone is enough to determine the validity of the deformed grid.
\end{abstract}


\maketitle

\section{Introduction}
The classical Lagrangian scheme solves the advection equation $f_t + {\bf u}\cdot \nabla f = 0$ by first tracing back the grid points and then estimating the value of $f$ at the traced-back grid points at the previous time step. The motivation behind it is that $f(t,{\bf q}(t))$ must be constant along every trajectory $(t,{\bf q}(t))$ of the flow particles if $f$ satisfies the advection equation. A natural generalization of the advection equation is the Lie-advection, or Lie transport: $\omega_t + \mathcal{L}_{\bf u}\omega = 0$, where $\omega(t,{\bf x})$ refers to a spatially and temporally varying differential form, and $\mathcal{L}_{\bf u}\omega$ reads as the Lie derivative of $\omega$ along the flow velocity field ${\bf u}$.  The conservation of quantity along trajectory lines is then generalized to conservation of quantities along trajectories of points, curves, faces, and bodies. Thus natural generalizations of the classical Lagrangian scheme include the geometric Lie transport/Lie-advection scheme \citep{Elcottetal2007}, the generalized Lagrangian scheme \cite{Hjort2017CSLAM,Lauritzen2018JAMES,Lauritzen2009JCP}, and some other methods that have this idea as part of their theoretical foundation \citep{Heumann2011,Heumann2013FCM,Heumann2012FullyDiscrete,Mao2026StructurePreserving,Tonnon2024SemiLagrangianFE} . First of all,  a differential k-form $\omega_k(t,{\bf x})$ is discretized as a collection of integrals over all k-dimensional elements of the grid $K$: $\omega_k^{\mathrm{h}}(t,{\bf x}) = \{\int_{\sigma_k}\omega_k(t): \sigma_k\in K\}$.  Then in order to estimate $\omega_k^{\mathrm{h}}(t+\Delta t)$, the generalized Lagrangian scheme performs the following two steps: 1) estimate all the time-backward deformed k-dimensional elements $\sigma_k^-$ for all $\sigma_k\in K$; 2) estimate $\int_{\sigma_k^-}\omega(t)$. The motivation behind the generalized Lagrangian scheme is similar to that for the classical Lagrangian scheme: $\int_{\sigma_k^-}\omega_k(t) = \int_{\sigma_k}\omega(t+\Delta t)$ should hold exactly as long as the flow is Lipschitz continuous over the time interval $[t,t+\Delta t]$, since a smooth flow map which preserves the boundary of the domain is always an isomorphism. At the discrete level, the deformation map $\phi(s,\cdot)$ is constructed manually. The artificial $\phi(\cdot,\cdot)$ is not a flow map in most of the cases. Thus a natural question is how to determine whether the artificial $\phi(s,\dot)$ is a global isomorphism. This question is also of concern in the community of arbitrary Lagrangian-Eulerian methods \citep{berndt2010using,breil2025mesh, cirrottola2021adaptive,dobrev2020simulation,guisset2024cell,johnen2013geometrical,kucharik2020improved,liu2023topology,toulorge2013robust,zhang2025interface}. More specifically, the arbitrary Lagangian-Eulerian (ALE) community needs a reliable and efficient method for determining whether, after an $\mathcal{O}(1)$ time interval, the deformed grid self-intersects or whether any element reverses its orientation. Since the discretized deformation is constructed by users in practice, it is reasonable to assume that a spatially and temporally piecewise-smooth deformation map $\phi: [0,\infty)\times \Omega\to\Omega$ is explicitly given, where $\Omega$ refers to the domain. If $\Omega$ has boundary, it is reasonable to assume further that $\phi$ maps $\partial\Omega$ to $\partial\Omega$. Then we have the following result:
\newtheorem{theorem}{Theorem}
\newtheorem{proposition}{Proposition}
\newtheorem{corollary}[proposition]{Corollary}
\newproof{pf}{Proof}
\newproof{pol1}{Proof of Lemma \ref{thm: lemma1}}
\newproof{pol2}{Proof of Lemma \ref{thm: lemma2}}
\newproof{poc}{Proof of Corollary \ref{thm: cor Prokhorova}}
\newproof{pot}{Proof of Theorem \ref{thm: criterion}}
\newtheorem{lemma}{Lemma}

\begin{theorem}\label{thm: criterion}
    Suppose $K$ is a grid for the compact domain $\Omega$. $\mathrm{dim}\Omega \leq 3$. For each closed grid element $K_i$ and for any $s \geq 0$, $\phi_{i}(s,\cdot): K_i\to\Omega$ is a smooth map. Assume that $\phi_i$ is piecewise-smooth with respect to $s$. $\phi_i$ and $\phi_j$ coincide on the common boundary of $K_i$ and $K_j$ for all $i,j$, i.e., $\phi_{i}\big{|}_{K_i\cap K_j} = \phi_j\big{|}_{K_i\cap K_j}$. Thus the global deformation map $\phi(s,\cdot):\Omega \to \Omega$ obtained by gluing the maps $\phi_i(s,\cdot)$ is well-defined. Assume further that  $s_0 > 0$ and the following conditions hold:
    \begin{enumerate}
    \itemsep=0pt
    \item $\phi(s,\partial \Omega)\subset \partial \Omega$ for every $s\in [0,s_0]$;
    \item The Jacobian determinant of $\phi_i(s_0,\cdot)$ is positive for every point in the closed element $K_i$ ;
    \item $\phi_i(0,\cdot) = \text{id$_{K_i}$}$ the identity map on $K_i$. 
    \end{enumerate}  
    Then $\phi(s_0,\cdot)$ is an isomorphism on $\Omega$. 
\end{theorem}

While it is obvious that positive Jacobian determinant in each grid element is a necessary condition for $\phi(s,\cdot)$ being an isomorphism, it is not in general true that local isomorphism would imply global isomorphism \cite{Garanzha2021, Weber2014}. Theorem \ref{thm: criterion} claims that, given the natural condition $\phi(s,\partial\Omega)\subset\partial\Omega$ for all $s$, the positivity of the Jacobian determinant is also a sufficient condition. To the best of our knowledge, any equivalent form of this result has not been explicitly stated in previous literature. However, Theorem \ref{thm: criterion} is closely related to a broader body of mathematical research on determining whether a local homeomorphism $f$ is also a global homeomorphism. For instance, it is proved in \cite{ lipman2014bijective, Massey1992, MeistersOlech1963, Xu2011, prokhorova2008, prokhorova2015} that $f$ is a global homeomorphism given that $f$ restricted to the boundary is injective or bijective, or has degree 1. However, these requirements on the boundary are theoretical and can not be verified directly. Theorem \ref{thm: criterion} essentially proves that, these boundary requirements can be eliminated in the grid-deformation setting. 


\section{The proof of Theorem \ref{thm: criterion}}
The proof presented here is based on the following result \citep{prokhorova2008,prokhorova2015}. 
\begin{theorem}[Prokhorova]
\label{thm: Prokhorova}
Suppose $X$ and $Y$ are connected topological manifolds of equal dimensions. X is compact. $f: X\to Y$ is a continuous mapping, $f(\partial X)\subset \partial Y$. $f\big{|}_{\partial X}$ and $f\big{|}_{X\backslash\partial X}$ are immersions. Then $f$ is a finite-fold covering. $f$ is an isomorphism (i.e., one-to-one and onto) if any of the following conditions is satisfied further:
\begin{enumerate}
\itemsep=0pt
\item there exists a point $y\in Y$ so that $f^{-1}(y)$ contains only one point;
\item $\partial X$ and $\partial Y$ are non-empty, and $f\big{|}_{\partial X}$ is injective;
\item (some other conditions that are more complicated and hard to verify).
\end{enumerate}  
\end{theorem}

The results of Prokhorova are used for grid-generation. Another work for grid-generation is \cite{lipman2014bijective}. In the grid-deformation setting, supposing that the deformation map is represented by $\phi$, we naturally have the boundary-preserving homotopy: $\phi(s,\cdot)\simeq_{\partial}\mathrm{id}_{\Omega}$, where $\simeq_{\partial}$ means that all the intermediate maps in the homotopy map $\partial \Omega$ to $\partial \Omega$.  This quickly leads to the following corollary of Prokhorova's theory.
\begin{corollary}\label{thm: cor Prokhorova}
    Suppose $K$ is a grid for a compact domain $\Omega$. For each closed grid element $K_i$ and for any $s \geq 0$, $\phi_{i}(s,\cdot): K_i\to\Omega$ is a smooth map. Assume that $\phi_i$ is piecewise-smooth with respect to $s$. $\phi_i$ and $\phi_j$ coincide on the common boundary of $K_i$ and $K_j$ for arbitrary $i,j$, i.e., $\phi_{i}\big{|}_{K_i\cap K_j} = \phi_j\big{|}_{K_i\cap K_j}$. Thus the global deformation map $\phi(s,\cdot):\Omega \to \Omega$ gluing together all $\phi_i(s,\cdot)$ is well-defined. Assume further that $s_0 > 0$ and the following conditions hold:
    \begin{enumerate}
    \itemsep=0pt
    \item $\phi(s,\partial \Omega)\subset \partial \Omega$ for any $s\in [0,s_0]$;
    \item $\phi(s_0,\cdot)$ is a local isomorphism at all points in $\Omega$;
    \item $\phi_i(0,\cdot) = \text{id$_{K_i}$}$ is the identity map.
\end{enumerate}  
    Then $\phi(s_0,\cdot)$ is a global isomorphism. 
\end{corollary}
\begin{poc}
    Given that $\phi(s_0,\cdot)$ is a local isomorphism at every point in $\Omega$ and that $\phi(s_0,\partial\Omega)\subset\partial \Omega$, Theorem \ref{thm: Prokhorova} implies that $\phi(s_0,\cdot)$ is a finite-fold covering. Since $\phi(s,\cdot)$ maps boundary to the boundary for all $s\in[0,s_0]$, $\phi(s_0,\cdot)\simeq_{\partial} \mathrm{id}_{\Omega}$. Thus $\mathrm{deg}(\phi(s_0,\cdot))=\mathrm{deg}(\mathrm{id}_{\Omega}) = 1$. Since $\phi(s_0,\cdot)$ is a local isomorphism, the local degree $\mathrm{deg}(\phi(s_0,\cdot),{\bf x})$ must be a constant on each connected component. For each ${\bf y}\in\Omega$, $\phi(s_0,\cdot)^{-1}({\bf y})$ consists of finitely many isolated points $\{{\bf x}_1,{\bf x}_2, \cdots, {\bf x}_m \}$ since $\Omega$ is compact. These points must lie in the same connected component since $\phi(s_0,\cdot)\simeq\mathrm{id}_{\Omega}$. Then the degree formula $\mathrm{deg}(\phi(s_0,\cdot)) = \displaystyle\sum_{i=1}^m\mathrm{deg}(\phi(s_0,\cdot),{\bf x}_i,{\bf y})$
    implies that $m=1$ and $\mathrm{deg}(\phi(s_0,\cdot),{\bf x}_1,{\bf y}) = 1$. Therefore $\phi(s_0,\cdot)$ is a global isomorphism. This completes the proof of Corollary \ref{thm: cor Prokhorova}.\qed
\end{poc} 

\begin{lemma}\label{thm: lemma1}
    Under the assumptions of Theorem \ref{thm: criterion}, we have the following results:
    \begin{enumerate}
    \itemsep=0pt
    \item $\phi^{-1}(\mathrm{y})$ are isolated points which are not on the boundary of any element for almost every $\mathrm{y}\in \Omega$;
    \item $|\phi^{-1}(\mathrm{y})| = 1$ for any such $\mathrm{y}$;
    \item for almost every $\mathrm{x}\in\Omega$, $|\phi(s_0,\cdot)^{-1}(\phi(s_0,\mathrm{x}))|=1$.
    \end{enumerate}
\end{lemma}
\begin{pol1}
    This is a direct consequence of Sard's embedding theorem \citep{Sard1942}. Since the boundary of the elements in $K$ form a lower dimensional piecewise-smooth manifold, the image of all the faces of $K$ under $\phi(s_0,\cdot)$ has Lebesgue measure zero. Thus the points in $\phi(s_0,\cdot)^{-1}(\mathrm{y})$ do not lie on any boundary of any element for  every $\mathrm{y}$ not in the image of the faces under $\phi(s_0,\cdot)$. For any such $\mathrm{y}$, we have the degree formula $\mathrm{deg}(\phi(s_0,\cdot)) = \sum_{\mathrm{x}\in\phi(s_0,\cdot)^{-1}(\mathrm{y})}\mathrm{deg}(\phi(s_0,\cdot),\mathrm{x},\mathrm{y})$. Since the Jacobian determinant is positive at these $\mathrm{x}$'s, $\mathrm{deg}(\phi(s_0,\cdot),\mathrm{x},\mathrm{y}) = 1$. Thus $\mathrm{deg}(\phi(s_0,\cdot)) = |\phi(s_0,\cdot)^{-1}(\mathrm{y})|$. Since $\phi(s_0,\cdot)\simeq_\partial \mathrm{id}_{\Omega}$, $\mathrm{deg}(\phi(s_0,\cdot)) = \mathrm{deg}(\mathrm{id}_{\omega}) = 1$. This suggests that $|\phi(s_0,\cdot)^{-1}(\mathrm{y})| = 1$. The third statement is a direct consequence of the previous two statements. \qed
\end{pol1}

\begin{flushleft}
    {\bf Assumption on the regularity of $K$:}We assume that $K$ is a regular grid for $\Omega$. More specifically, any two distinct faces are not tangent at any common point; any two distinct edges are not tangent at any common point; a face and an edge are not tangent unless the edge is contained in the face.
\end{flushleft}
\begin{pot}
     We only need to prove that $\phi(s_0,\cdot)$ is a local isomorphism. Then Corollary \ref{thm: cor Prokhorova} would imply Theorem \ref{thm: criterion}. Since the Jacobian determinant of $\phi(s_0,\cdot)$ is positive at the interior points of each element, $\phi(s_0,\cdot)$ is already a local isomorphism at these points. Thus all we need is to prove that $\phi(s_0,\cdot)$ is a local isomorphism on the boundary points of each closed element. We only prove this result for the case $\mathrm{dim}\Omega = 3$. The proof for lower-dimensional cases is similar and easier, and is thus left to the reader. The proof is performed case-by-case as follows.  

    \begin{flushleft}
        {\bf Case 1} (${\bf x}$ is an interior point at a face of some element):
    \end{flushleft}

    If ${\bf x}\in\partial \Omega$, then there is nothing to prove. Now assume that ${\bf x}$ is an interior point of $\Omega$. This means that  ${\bf x}$ belongs to a common face $\sigma$ of two elements, denoted by $K_i$ and $K_j$, but ${\bf x}$ is not on any edge of $K_i$ or $K_j$. Let $U$ be a small ball centered at ${\bf x}$. Denote by $U_i := U\cap K_i$, $U_j := U\cap K_j$. Then $U_i\cong\phi_i(s_0,U_i)$, $U_j\cong\phi_j(s_0,U_j)$, $U_i\cap \partial K_i\cong\phi_i(s_0,U_i\cap\partial K_i) = \phi_j(s_0,U_j\cap\partial K_j) \cong U_j\cap\partial K_j$. Denote  by $\sigma_U^{\phi} = \phi_i(s_0,U_i\cap\partial K_i)$. Since $\phi_i(s_0,\cdot)$ is smooth on the closure of $K_i$, $\sigma^\phi_U$ is a smooth surface in $\Omega$ with well-defined orientation. And for each point ${\bf y}$ in $\phi(s_0,U)\backslash\sigma_U^{\phi}$, it can be well-defined which side of $\sigma_U^\phi$ ${\bf y}$ lies in. Since the Jacobian determinant of $\phi_i(s_0,\cdot)$ and $\phi_j(s_0,\cdot)$ are both positive, $\phi_i(s_0, U_i)$ and $\phi_j(s_0,U_j)$ must be at different sides of $\sigma_U^\phi$. Therefore $\phi_i(s_0,U_i)\cap\phi_j(s_0,U_j) = \sigma^\phi_U$. Thus $\phi(s_0,\cdot)$ is a local isomorphism at ${\bf x}$.
    \begin{flushleft}
        {\bf Case 2} (${\bf x}$ is an interior point at an edge of some element):        
    \end{flushleft}

    This means that ${\bf x}\in e$ for some edge $e\in K$, but ${\bf x}$ is not a vertex of $K$. If ${\bf x}\in\partial \Omega$, then ${\bf x}$ is on an edge of a two-dimensional grid for the two-dimensional domain $\partial \Omega$. This case is essentially the same as the case 1 considered above but for the lower-dimensional domain $\partial\Omega$. The proof is similar, and thus left to the reader. In what follows, it is assumed that ${\bf x}$ is an interior point of the three-dimensional domain $\Omega$. ${\bf x}$ could belong to several different elements at the same time. Denote these elements by $K_0,...,K_{L}$, so that $K_i$ and $K_{\mathrm{mod}(i+1,L+1)}$ share a common face containing the edge $e$.  For each of these elements $K_i$, $K_i$ has two faces that contain $e$. Denote all these faces by $\sigma_0,\sigma_1,\cdots,\sigma_{L}$. The indices of these faces are chosen so that $\sigma_i,\sigma_{\mathrm{mod}(i+1,L+1)}\in K_i$. In this case, all indices in the following discussion are taken modulo $L+1$. Equivalently, $i+1$ is understood as $(i+1)\bmod(L+1)$ whenever $i+1>L$. Let $U$ be a small open tube containing ${\bf x}$ and whose central axis is part of $e$. Denote by $U_i := U\cap K_i$,  $\sigma_{U,i}:= \sigma_i\cap U$, $e_U:= e\cap U$, $\sigma^\phi_{U,i} := \phi(s_0,\sigma_{U,i})$, $e^\phi_U := \phi(s_0,e_U)$. Choose a local orientation-preserving coordinate $\Phi_{1}$ for ${\bf x}$  satisfying the following conditions:
    \begin{enumerate}
    \itemsep=0pt
    \item $\Phi_{1}: (-\eta,\eta)\times D_{\epsilon}\to U$ is a diffeomorphism, where $D_\epsilon$ refers to the two-dimensional disk of radius $\epsilon$; 
    \item $\Phi_1((0,(0,0))) = {\bf x}$;\hspace{2mm} $\Phi_1((-\eta,\eta)\times \{(0,0)\}) = e_U$;
    \item $\Phi_1$ restricted to $((-\eta,\eta)\times \{(0,0)\})$ is orientation-preserving;
    \item For each $i$, $\Phi_{1}^{-1}(\sigma_{U,i}) = (-\eta,\eta)\times \{(r,\theta_i)\in D_\epsilon: r\in [0,\epsilon]\}$ for some $\theta_i\in [0,2\pi)$
    \end{enumerate}  
    By properly choosing the order of the indices $\{0,1,\cdots,L\}$, we may assume that $0=\theta_0 < \theta_1 < \cdots < \theta_L < 2\pi$. Since $U$ is small, $K$ is regular, and $\phi_i(s,\cdot)$ is smooth, $\sigma_{U,i}^\phi$ has well-defined normal vectors and the oriented normal directions of $\sigma_{U,i}^\phi$ and $\sigma_{U,j}^\phi$ are well-separated. Let ${\bf y}_{s_0}= \phi(s_0,{\bf x})$. Choose a local orientation-preserving coordinate $\Phi_{s_0,2}$ for ${\bf y}_{s_0}$ so that the following conditions are satisfied:
    \begin{enumerate}
    \itemsep=0pt
    \item $\Phi_{s_0,2}:(-\eta,\eta)\times D_{\epsilon}\to \Omega$ is a smooth embedding;
    \item $\Phi_{s_0,2}((0,(0,0))) = {\bf y}_{s_0}$;\hspace{2mm} $\Phi_{s_0,2}((-\eta,\eta)\times \{(0,0)\}) = e_U^\phi$;
    \item For each $i$, $\Phi_{s_0,2}^{-1}(\sigma_{U,i}^\phi) \subset (-\eta,\eta)\times \{(r,\theta_i^\phi)\in D_\epsilon: r\in [0,\epsilon]\}$ for some $\theta_i^\phi\in [0,2\pi)$.
    \end{enumerate}  
    We point out that such a local coordinate $\Phi_{s_0,2}$ exists even if $\phi(s_0,\mathrm{x})\in\partial \Omega$ is logically allowed. 
    Choose small $\eta_{1}$ and $\epsilon_{1}$ so that $\phi(s,\Phi_1((-\eta_1,\eta_1)\times D_{\epsilon_1})) \subset \Phi_{s,2}((-\eta,\eta)\times D_{\epsilon})$. Then we have the following piecewise-smooth map:
    \begin{align}
        \psi:=\Phi_{s_0,2}^{-1}\circ\phi(s_0,\cdot)\circ\Phi_{1}:(-\eta_1,\eta_1)\times D_{\epsilon_1}\to (-\eta,\eta)\times D_{\epsilon}.\nonumber
    \end{align}
    It is assumed that $\Phi_{s_0,2}$ is carefully chosen so that $\psi$ restricted to $(-\eta_1,\eta_1)\times (0,0)$ is increasing. 
    Let $\psi_i$ be the restriction of $\psi$ to $(-\eta_1,\eta_1)\times \{(r,\theta)\in D_{\epsilon_1}: r\in [0,\epsilon_1], \theta\in[\theta_i,\theta_{i+1}]\}$. We have the following results:
    \begin{enumerate}
    \itemsep=0pt
    \item $\psi_i$ is a smooth orientation-preserving embedding;
    \item $\psi_i$ maps $(-\eta_1,\eta_1)\times \{(r,\theta_i)\in D_{\epsilon_1}: r\in [0,\epsilon_1]\}$ into $(-\eta,\eta)\times \{(r,\theta_i^\phi)\in D_{\epsilon}: r\in [0,\epsilon]\} $
    \end{enumerate}  
    Let $\theta$ be a variable that changes from $0$ to $2\pi$. Then for each $a\in (-\eta_1,\eta_1)$ and $r\in (0,\epsilon_1)$. Denote by $\gamma_{a,r,1} = \{(a,(r,\theta)):\theta\in [0,2\pi]\}$, which is a closed curve in $(-\eta_1,\eta_1)\times D_{\epsilon_1}$, and denote by $\gamma_{a,r,2} := \{(a,(r,\theta)): \theta\in [0,2\pi]\}$ the image of $\gamma_{a,r,1}$. According to case 1, $\psi$ restricted to $\gamma_{a,r,1}$, denoted by $\psi_{a,r}$, is a local isomorphism. Theorem \ref{thm: Prokhorova} suggests that $\psi_{a,r}$ is then a covering map. The covering degree of $\psi_{a,r}$ must the the same integer, denoted by $n(\mathrm{x})$, for small enough $a$ and $r$.  Since $\psi_{a,r}$ is an orientation-preserving local isomorphism, the local degree of $\psi_{a,r}$ must  equal to one. If $n(\mathrm{x})\neq 1$, then  $|\psi_{a,r}^{-1}(\mathrm{y})| > 1$ for any $\mathrm{y}\in \gamma_{a,r,2}$. This means that $|\phi(s_0,\cdot)^{-1}(\psi(s_0,\mathrm{x}))| > 1$ for all $\mathrm{x}\in U\backslash e$, which contradicts with Lemma \ref{thm: lemma1}. Thus $n(\mathrm{x}) = 1$ and $\psi$ restricted to $(-\eta_1,\eta_1)\times D_{\epsilon_1}\backslash\{(0,0)\}$ is injective. This further induces that $\psi$ is injective on the whole $(-\eta_1,\eta_1)\times D_{\epsilon_1}$. Hence $\phi(s_0,\cdot)$ restricted to $U$ is injective.  The theorem of invariance of domain \citep{brouwer1912invariance} implies that $\phi(s_0,\cdot)$ is an isomorphism on $U$. This completes the proof of case 2.

    \begin{flushleft}
        {\bf Case 3} (${\bf x}$ is a vertex of some element):
    \end{flushleft}

    If ${\bf x}\in\partial \Omega$, the idea of the proof is similar to, but much simpler than the proof of case 2. Thus we leave it to the interested readers. Now assume that ${\bf x}$ is an interior point of $\Omega$. Let $U$ be a small ball centered at ${\bf x}$. The point ${\bf x}$ belongs to several distinct elements, which we denote by $K_0, K_1,\cdots, K_L$. Since $U$ is small, ${\bf x}$ is the only vertex in $U$. Case 1 and case 2 suggest that $\phi(s,\cdot)$ is already a local isomorphism on $U_0 := U\backslash\{{\bf x}\}$. Let ${\bf y}_{s_0} = \phi(s_0,{\bf x})$. Logically, it is allowed that $\mathrm{y}_{s_0}\in\partial\Omega$. Then ${\bf x}$ must be an isolated point of $\phi(s_0,\cdot)^{-1}(\mathrm{y}_{s_0})$ since only finitely many elements contain $\mathrm{x}$ and  $\phi(s_0,\cdot)$ restricted to each $K_i$ is assumed to be a local isomorphism. For any small ball $V\subset U$ centering at $\mathrm{x}$, $\phi(s_0,\cdot)$ restricted to $\partial V$ is then a local isomorphism. Similar to the analysis presented in case 2 and apply Theorem \ref{thm: Prokhorova} again, it can be proved that $\phi(s_0,\cdot)$ restricted to $\partial V$ is an isomorphism for any such $V$. Thus $\phi(s_0,\cdot)$ restricted to $U$ is injective. The theorem of invariance of domain \citep{brouwer1912invariance} implies that $\phi(s,\cdot)$ is an isomorphism on $U$. This completes the proof of case 3, hence of the full theorem. \qed
\end{pot}
\section{Summary and discussion}
Corollary \ref{thm: cor Prokhorova} provides a rigorous theoretical criterion that follows directly from Prokhorova's theory and the natural homotopy in grid deformation. However, the second condition in Corollary \ref{thm: cor Prokhorova} is usually not easy to check directly on the boundary of each element $K_i$, since the deformation map constructed manually is usually piecewise-smooth. In fact, without the natural homotopy in grid deformation, being an isomorphism on each element would not imply being a local isomorphism everywhere. This practical difficulty is eliminated in Theorem \ref{thm: criterion}. The main contribution of this manuscript is to show that, in the grid-deformation setting, a local isomorphism restricted to each element automatically induces a local isomorphism on the boundary of each element. Theorem \ref{thm: criterion} also applies to $\Omega$ without boundary, a case that cannot be handled by the methods in \cite{lipman2014bijective}.

Theorem \ref{thm: criterion} shows that $\phi(s_0,\cdot)$ is an isomorphism under certain conditions. However, Theorem \ref{thm: criterion} does not guarantee that the deformed grid for $s\in (0,s_0)$ is also valid. In principle, it is possible that $\phi(s,\cdot)$ becomes degenerate at some time in $(0,s_0)$ and then recovers at time $s_0$. Thus, an interesting question is: under the conditions of Theorem \ref{thm: criterion}, does there exist another deformation process $\tilde{\phi}(s,\cdot)$ such that 1) $\phi$ and $\tilde{\phi}$ are homotopic; 2) $\phi(s_0,\cdot) = \tilde{\phi}(s_0,\cdot)$; 3) $\tilde{\phi}(s,\cdot)$ is always an isomorphism with positive Jacobian determinant in each closed element for every $s\in [0,s_0]$? The answer to this question is affirmative for two-dimensional domains \cite{Epstein1966} and negative in general \cite{FriedmanWitt1986}. We leave the study of this question in practical scenarios for future work.

\section*{Acknowledgements}
Funding: this work is supported by the National Natural Science Foundation of China [grant numbers 42350003, 12671486].
\printcredits

\bibliographystyle{cas-model2-names}

\bibliography{Mendeley}

\end{document}